\documentclass{amsart}
\usepackage{amsfonts}
\usepackage{amsmath,amssymb}
\usepackage{amsthm}
\usepackage{amscd}
\usepackage{graphics}
\usepackage{graphicx}

\theoremstyle{remark}{
\newtheorem{Def}{{\rm Definition}}
\newtheorem{Ex}{{\rm Example}}
\newtheorem{Rem}{{\rm Remark}}

}
\theoremstyle{plain}{

\newtheorem{Prop}{Proposition}

\newtheorem*{MainThm}{Main Theorem}
\newtheorem*{MainCor}{Main Corollary}

}

\begin{document}
\title[An explicit differential topological study of compact manifolds]{A differential topological study of compact manifolds having simple structures}
\author{Naoki Kitazawa}
\address{Institute of Mathematics for Industry, Kyushu University, 744 Motooka, Nishi-ku Fukuoka 819-0395, Japan\\
 TEL (Office): +81-92-802-4402 \\
 FAX (Office): +81-92-802-4405 \\
}
\email{n-kitazawa@imi.kyushu-u.ac.jp}
\urladdr{https://naokikitazawa.github.io/NaokiKitazawa.html}
\maketitle
\begin{abstract}
The present paper mainly presents, for example, explicit classifications of compact smooth manifolds having non-empty boundaries and simple structures where the dimensions are general. 
Studies of this type is fundamental and important. They also remain to be immature and difficult. This is due to the fact that the dimensions are high and this has prevented us from studying the manifolds in geometric and constructive ways.

Moreover, most of the present work is motivated by explicit studies of higher dimensional variants of Morse functions: especially so-called {\it special generic} maps.
The class of special generic maps is a natural class containing canonical projections of unit spheres and Morse functions on homotopy spheres with exactly two singular points. Their images are in general (compact) manifolds smoothly immersed to the targets and the dimensions of the images and the targets coincide. 
They know much about the topologies and the differentiable structures of the manifolds of the domains. The author has previously studied related problems and the present study gives new extended results of related previous results. Last, we also present a dream for contribution to studies of special generic maps and higher dimensional variants of Morse functions and manifolds admitting them.

\end{abstract}


\maketitle
\section{Introduction.}
\label{sec:1}
Compact manifolds which may have non-empty boundaries are fundamental geometric objects. The class is more and more general than the class of closed manifolds. It has been also difficult to classify explicit classes of closed manifolds. 
$1$-dimensional and $2$-dimensional cases are classical. In the last century, closed and simply-connected manifolds whose dimensions are greater than $4$ are classified mainly via sophisticated algebraic topological methods and differential topological ones such as homotopy theory, surgery theory, Morse theory, and so on: \cite{milnor2}, \cite{milnor3}, \cite{wall2}, and so on, explain related classical theory for example. $3$ or $4$-dimensional closed manifolds are still difficult and interesting. It is also difficult to understand closed and simply-connected manifolds whose dimensions are greater than $4$ in geometric and constructive ways due to the assumption that the dimensions are high and Nishioka and the author has obtained related results, mainly for $5$ or $7$-dimensional closed and simply-connected manifolds, via higher dimensional variants of Morse functions with exactly two singular points on {\it homotopy spheres} and more general smooth maps: a {\it homotopy sphere} means a smooth manifold homeomorphic to a unit sphere in the present paper. Nishioka's work is \cite{nishioka} and related works of the author are some of \cite{kitazawa0.1}--\cite{kitazawa13}. 
Especially, Nishioka's result is on $5$-dimensional closed and simply-connected manifolds and such smooth maps on them and closely related to the main ingredient of the present paper. Last, if we remove the restriction that the manifolds are simply-connected, then it seems to be far more difficult even for algebraic or abstract studies.

In the present paper, we study compact, simply-connected and smooth manifolds of suitable classes with non-empty boundaries in geometric and constructive ways. We present Main Theorem and Main Corollary with several undefined notions and notation, which are explained in the third section. 

\begin{MainThm}
Let $n>k>1$ and $a>0$ be integers.
\begin{enumerate}
\item
\label{m:1}
Let $K$ be a ${\rm (}k-1{\rm )}$-connected polyhedron whose dimension is smaller than or equal to $n-a$.  
Assume that the values of a root $s_0$ for $K$ are always of the form $({[Y]}_{{\rm Diff}},1,p)$ satisfying either of the following two.
\begin{enumerate}
\item $({[Y]}_{{\rm Diff}},1,1)$ for a closed and connected manifold $Y$ which can be immersed {\rm (}resp. embedded{\rm )} smoothly into ${\mathbb{R}}^{\dim Y+a}$ with a trivial normal bundle.
\item $({[Y]}_{{\rm Diff}},1,0)$ for a closed and connected manifold $Y$ which can be immersed {\rm (}resp. embedded{\rm )} smoothly into ${\mathbb{R}}^{n}$.
\end{enumerate}

Then $Y$ here is always {\rm (}$k-1${\rm )}-connected, there exists an $n$-dimensional compact, connected and smooth manifold $X$ and $X$ is an {\rm SIE-$K$} {\rm (}resp. {\rm SEE-$K$}{\rm )}.
\item
\label{m:2}
For an $n$-dimensional compact, connected and smooth manifold $X$ and an elementary polyhedron $K$, suppose that $X$ is an SIE-$K$ {\rm (}SEE-$K${\rm )} and that the following two are satisfied.
\begin{enumerate}
\item $X$ and $K$ are {\rm (}$k-1${\rm )}-connected.
\item For a root $s_0$ for $K$, values are always of either of the following forms.
\begin{enumerate}
\item $({[Y]}_{{\rm Diff}},1,1)$ for a closed and connected manifold $Y$ which can be immersed {\rm (}resp. embedded{\rm )} smoothly into ${\mathbb{R}}^{\dim Y+a}$ with a trivial normal bundle.
\item $({[Y]}_{{\rm Diff}},1,0)$ for a closed and connected manifold $Y$ which can be immersed {\rm (}resp. embedded{\rm )} smoothly into ${\mathbb{R}}^{n}$.
\end{enumerate}
\end{enumerate}
Then $K$ is represented as a polyhedron obtained by a finite iteration of taking a suitable bouquet starting from polyhedron satisfying at least one of the following two.
\begin{enumerate}
\item A product of the following two polyhedra.
\begin{enumerate}
\item A closed and {\rm (}$k-1${\rm )}-connected manifold $F$ we can smoothly immerse {\rm (}resp. embed{\rm )} into ${\mathbb{R}}^{\dim F+a}$ with a trivial normal bundle.
\item A polyhedron $K_F$ such that for a root $s_{0,F}$ for $K_F$, values are always of the form $({[Y]}_{{\rm Diff}},1,1)$ for a closed and connected manifold $Y$ which can be immersed {\rm (}resp. embedded{\rm )} smoothly into ${\mathbb{R}}^{\dim Y+a}$ with a trivial normal bundle and that cannot be represented as a product of a closed and {\rm (}$k-1${\rm )}-connected manifold $Y^{\prime}$ we can smoothly immerse into ${\mathbb{R}}^{\dim Y^{\prime}+a}$ and another polyhedron $K_{Y,F}$.
\end{enumerate}
\item A closed and {\rm (}$k-1${\rm )}-connected manifold $F_0$ we can smoothly immerse {\rm (}resp. embed{\rm )} into ${\mathbb{R}}^{n}$.
\end{enumerate}
\end{enumerate}
\end{MainThm}
\cite{kitazawa13} concerns the case the following two are satisfied.
\begin{enumerate}
\item $n \leq 3k$.
\item The value of a root is always of the form $({[Y]}_{{\rm Diff}},1,p)$ where $Y$ is a homotopy sphere except for manifolds playing same roles as the manifold ''$F_0$'' plays here.
\end{enumerate}
\begin{MainCor}
In Main Theorem, let $n=7$ and $k=2$. Then in {\rm (}\ref{m:2}{\rm )}, $F$, $F_0$ and $K_F$ satisfy the following properties.
\begin{enumerate}
\item $F$ is a $k_F$-dimensional closed and simply-connected manifold with $2 \leq k_F \leq 4$.
\item $K_F$ is represented as a bouquet of copies of unit spheres if $k_F \geq 3$ and in general represented as a bouquet of closed and simply-connected manifolds whose dimensions are smaller than $7-k_F$.
\item $F_0$ is represented as a connected sum of finitely many closed and simply-connected manifolds whose dimensions are smaller than $7$ we can smoothly immerse {\rm (}resp. embed{\rm )} into ${\mathbb{R}}^n$. 
\end{enumerate}
\end{MainCor}

\cite{kitazawa13} concerns the case $n=5,6$.

The present study is also motivated by studies of higher dimensional variants of Morse functions, especially ones on {\it special generic} maps.

Hereafter, a {\it singular} point $p \in X$ of a differentiable map $c:X \rightarrow Y$ is a point at which the differential ${dc}_p$ satisfies ${\rm rank}{dc}_p<\min\{\dim X,\dim Y\}$. Let $S(c)$ denote the set of all singular points of $c$ and we call this the {\it singular set} of $c$. We call $c(S(c))$ the {\it singular value set} of $c$ and the complementary set $Y-c(S(c))$ the {\it regular value set} of $c$. A {\it singular value} of $c$ is a point in the singular value set and a {\it regular value} of $c$ is a point in the regular value set.

A {\it special generic} map on an $m$-dimensional closed and smooth manifold into ${\mathbb{R}}^n$  is a smooth map at each singular point of which it is represented as $(x_1, \cdots, x_m) \mapsto (x_1,\cdots,x_{n-1},\sum_{k=n}^{m}{x_k}^2)$
for suitable coordinates. The restriction to the singular set, which is shown to be an ($n-1$)-dimensional closed and smooth submanifold with no boundary, is a smooth immersion. 
The simplest special generic maps are Morse functions with exactly two singular points on homotopy spheres, playing important roles in so-called Reeb's theorem, and the canonical projection of the unit sphere $S^m$ into ${\mathbb{R}}^n$ for $m\geq n \geq 1$ where $S^m$ denotes the $m$-dimensional unit sphere. The image of a special generic map is the image of a smooth immersion of an $n$-dimensional compact and smooth manifold into ${\mathbb{R}}^n$ and the $n$-dimensional manifold inherits important information of the closed manifold such as homology groups, cohomology rings, and so on, in considerable cases. 
Special generic maps restrict the topologies and the differentiable structures of the closed manifolds strongly in considerable cases. \cite{saeki} and \cite{saeki2} are pioneering studies and related studies of the author are \cite{kitazawa8} and \cite{kitazawa10} for example.

The author previously studied the topologies and the differentiable structures of the images of in \cite{kitazawa13} and also in \cite{kitazawa11} and \cite{kitazawa12} before for example. 
\cite{nishioka} motivated the author to study this. He studied the existence problem on special generic maps on $5$-dimensional closed and simply-connected manifolds, completely classified in the topology, PL, and smooth categories, in \cite{barden} and studied in \cite{smale}, and solved this completely. Through the study, he has shown a proposition stating that integral homology groups of $4$-dimensional compact, simply-connected and smooth (PL) manifolds with non-empty boundaries are free and that this is generalized to $k$-th integral homology groups of ($k+2$)-dimensional compact, simply-connected, and smooth manifolds with non-empty boundaries. As a main theorem, a $5$-dimensional closed and simply-connected manifold admits a special generic map into ${\mathbb{R}}^n$ if and only if it is homeomorphic and as a result diffeomorphic to $S^5$ or represented as a connected sum of total spaces of smooth bundles over the $2$-dimensional unit sphere $S^2$ whose fibers are diffeomorphic to $S^3$ for $n=3,4$ (such a manifold admits a special generic map into ${\mathbb{R}}^n$ if and only if it is homeomorphic to and as a result diffeomorphic to $S^5$ for $n=1,2$). This result implies that a $5$-dimensional closed and simply-connected manifold admitting a special generic map into ${\mathbb{R}}^n$ admits a special generic map into ${\mathbb{R}}^n$ for $n=1,2,3,4$ satisfying the following properties.
\begin{enumerate}
\item The restriction to the singular set is an embedding.
\item The image is diffeomorphic to the unit disc $D^n$ or diffeomorphic to a manifold represented as a boundary connected sum of finitely many copies of $S^2 \times D^{n-2}$ where $D^{k}$ denotes the $k$-dimensional unit disc. 
\end{enumerate}
 
The author has studied these images for higher dimensional cases in cited articles. Main Theorem and Main Corollary are extended new results to existing closely related results.

As a remark, in the introduction of \cite{kitazawa13}, importance of compact (and smooth) manifolds with non-empty boundaries, in the knot theory, low-dimensional geometry, and several explicit studies of algebraic geometry such as topological properties of algebraic curves, complementary spaces of affine subspaces of real or complex vector spaces, and so on, is also explained shortly. It is also said that general algebraic topological or differential topological theory of compact manifolds possibly with non-empty boundaries is immature.

In the next section, we add important properties of special generic maps. The third section is devoted to notions on compact manifolds with non-empty boundaries appearing as the domains of the immersions presented in the explanation on the images of special generic maps before and polyhedra they collapse to. They are essentially introduced first in \cite{kitazawa13}. We introduce these notions in a bit revised form. We prove Main Theorem in the fourth section and the fifth section is for a conclusion remark, saying that from the images and the structures of the special generic maps, we can know the manifolds of the domains. 

In the remaining sections, manifolds maps between manifolds, (boundary) connected sums, and so on, are smooth (of the class $C^{\infty}$), unless otherwise stated. 

\section{Special generic maps.}
\begin{Prop}[\cite{saeki}.]
\label{prop:1}
Let $m>n \geq 1$ be integers. 
An $m$-dimensional closed manifold $M$ admits a special generic map into ${\mathbb{R}}^n$ if and only if the following conditions hold.
\begin{enumerate}
\item There exists a 
smooth surjection $q_f:M \rightarrow W_f$ onto an $n$-dimensional compact manifold such that $q_f(S(f))= \partial W_f$ and an immersion $\bar{f}:W_f \rightarrow {\mathbb{R}}^n$ such that $f=\bar{f} \circ q_f$.
\item There exists a small collar neighborhood $C(\partial W_f)$ such that the composition of $q_f {\mid}_{{q_f}^{-1}(C(\partial W_f))}$ with the canonical projection to $\partial W_f$ gives a linear bundle whose fiber is diffeomorphic to $D^{m-n+1}$.
\item $q_f {\mid}_{{q_f}^{-1}(W_f-{\rm Int}(C(\partial W_f)))}$ gives a smooth bundle whose fiber is diffeomorphic to $S^{m-n}$.
\end{enumerate}

Furthermore, we can obtain a desired special generic map $f$ as $\bar{f} \circ q_f$ in this scene.
\end{Prop}
\begin{Ex}
\label{ex:1}
Let $m>n \geq 2$ and $l>0$ be integers. Let $\{S^{l_j} \times S^{m-l_j}\}_{j=1}^{l}$ be a family of products of unit spheres satisfying $1 \leq l_j \leq n-1$. Let  $M$ be a manifold represented as a connected sum of $l$ manifolds each of which is diffeomorphic to the corresponding product in the family $\{S^{l_j} \times S^{m-l_j}\}_{j=1}^{l}$. This admits a special generic map into ${\mathbb{R}}^n$ such that $\bar{f}$ is an embedding of a manifold represented as a boundary connected sum of the $l$ manifolds each of which is diffeomorphic to $S^{l_j} \times D^{n-l_j}$.
\end{Ex}
\begin{Prop}[\cite{nishioka}, \cite{saeki}, and so on.]
\label{prop:2}
\begin{enumerate}
\item Let $k>3$ be an integer. For a {\rm (}$k${\rm )}-dimensional compact and simply-connected manifold $X$ satisfying $\partial X \neq \emptyset$, $H_{k-2}(X;\mathbb{Z})$ is free. 
\item An $m$-dimensional closed and simply-connected manifold $M$ admits a special generic map $f$ into ${\mathbb{R}}^2$ for $m \geq 2$ if and only if it is homeomorphic to $S^m$ for $m \neq 4$ and diffeomorphic to $S^m$ for $m=4$. Furthermore, in this case, $M$ admits a special generic map $f$ into ${\mathbb{R}}^2$ such that $\bar{f}$ in Proposition \ref{prop:1} is an embedding of a copy of $D^2$ for $m \geq 3$.
\item If an $m$-dimensional closed and simply-connected manifold $M$ admits a special generic map $f$ into ${\mathbb{R}}^3$ for $m \geq 4$, then it admits a special generic map $f$ into ${\mathbb{R}}^3$ such that $\bar{f}$ in Proposition \ref{prop:1} is an embedding of a copy of $D^3$ or a manifold diffeomorphic to one represented as a boundary connected sum of finitely many copies of $S^2 \times D^1$.
\item A $5$-dimensional closed and simply-connected manifold $M$ admits a special generic map into ${\mathbb{R}}^n$ if and only if the following two hold. 
\begin{enumerate}
\item $M$ is homeomorphic to $S^5$ and as a result diffeomorphic to $S^5$ for $n=1,2$.
\item $M$ is homeomorphic to $S^5$ and as a result diffeomorphic to $S^5$, or homeomorphic to a manifold represented as a connected sum of total spaces of smooth bundles over $S^2$ whose fibers are diffeomorphic to $S^3$ and as a result diffeomorphic to this for $n=3,4$.
\end{enumerate}
Furthermore, in the latter case $n=3,4$, $M$ admits a special generic map $f$ into ${\mathbb{R}}^n$ such that $\bar{f}$ in Proposition \ref{prop:1} is as in Example \ref{ex:1} with $l_j=2$.
\end{enumerate}
\end{Prop}
\section{Compact manifolds with non-empty boundaries respecting $W_f$ and the image of a special generic map $f$ in Proposition \ref{prop:1} and important polyhedra.}
We can apply Proposition \ref{prop:1} to construct special generic maps from compact manifolds with non-empty boundaries smoothly immersed into the Euclidean spaces whose dimensions agree with the dimensions of the compact manifolds. 
Thus compact manifolds with non-empty boundaries respecting $W_f$ and the image of a special generic map $f$ in Proposition \ref{prop:1} are important objects. It is important to study their topologies and differentiable structures.

Notions and arguments in this section generalize several notions and arguments appearing first in \cite{kitazawa13}.
We call equivalence classes under the natural equivalence relation in the PL and the smooth category {\it PL types} and {\it diffeomorphism types}, respectively. For an object $X$ in the PL (smooth) category, we also call the equivalence class ${[X]}_{{\rm PL}} \ni X$ (resp. ${[X]}_{{\rm Diff}} \ni X$) the {\it PL {\rm (}resp. diffeomorphism{\rm )} type} for $X$. It is well-known that a smooth manifold $X$ is regarded as a polyhedron the PL type for which is canonically and uniquely defined. The PL type for $X$ means this. 

Let ${\mathbb{N}}_r:=\{x \in \mathbb{N} \mid x \leq r.\}$ for $r \in \mathbb{R}$. Let $s_0$ be a map on $\{0\} \sqcup {\mathbb{N}}_l$
or a sequence of lehgth $l+1$ indexed by the finite set for $l \in \mathbb{N}$ such that the values are of the form $({[Y]}_{{\rm Diff}},1,p)$: $Y$ is a smooth manifold, ${[Y]}_{{\rm Diff}}$ is the diffeomorphism type for $Y$, and $p=0,1$.
We consider an iteration of the following four procedures starting from $s_k:=s_0$.
\begin{enumerate}
\item
\label{s;1}
Take $s_k$ ($k \in \{0\} \sqcup {\mathbb{N}}_l$). 
\item
\label{s;2}
If $k=l$, then let $s_{l,0}:=s_l(0)$ and finish these four procedures.
\item
\label{s:3}
Choose two distinct numbers $k_1,k_2 \in \{0\} \sqcup {\mathbb{N}}_{l-k}$ with $k_1<k_2$. Define $s_{k+1}$ to be a map on $\{0\} \sqcup {\mathbb{N}}_{l-k-1}$
satisfying the following three conditions.
\begin{enumerate}
\item
\label{s:3a}
The values are of either of the following two forms.
\begin{enumerate}
\item $({[Y]}_{{\rm Diff}},1,p)$ for a smooth manifold $Y$ and $p=0,1$ where the first component is the diffeomorphism type for $Y$. 
\item $({[Y]}_{{\rm PL}},0,p)$ for a polyhedron $Y$ and $p=0,1$ where the first component is the PL type for $Y$.
\end{enumerate}
\item
\label{s:3b}
$s_{k+1}(j)=\begin{cases}
s_k(j) \quad {\rm(}0 \leq j \leq k_1-1{\rm )}\\
s_k(j+1) \quad {\rm(}k_1-1 < j \leq k_2-2{\rm )} \\
s_k(j+2) \quad {\rm(}k_2-2 < j < l-k-1{\rm )} 
\end{cases}.$
\item
\label{s:3c} $s_{k+1}(l-k-1)$ is a triplet of elements and satisfies either of the following three, including conditions on the original two elements $s_k(k_1)$ and $s_k(k_2)$.
\begin{enumerate}
\item
\label{s:3c1}
The first component is the PL type for a bouquet of a polyhedron in 
the first component of $s_k(k_1)$ and one in the first component of $s_k(k_2)$. The second component is the integer $0$. The third component is the product of the third component of $s_k(k_1)$ and that of $s_k(k_2)$.
\item
\label{s:3c2}
The first component is the PL type for a product of a smooth manifold and a polyhedron. 
The second component is $0$. The third component is $1$.
The polyhedra are a polyhedron or a smooth manifold in the first component of $s_k(k_1)$ and one in the first component of $s_k(k_2)$, respectively. For at least one of $s_k(k_1)$ and $s_k(k_2)$, the first component is a diffeomorphism type for a smooth manifold and the second component is $1$. The third component of $s_k(k_1)$ and that of $s_k(k_2)$ are both $1$.
\item
\label{s:3b3}
The first component is the diffeomorphism type for a smooth manifold represented as a connected sum of a smooth manifold in the first component of $s_k(k_1)$ and one in the first component of $s_k(k_2)$. The second component is $1$. The third component is the product of the third component of $s_k(k_1)$ and that of $s_k(k_2)$.
For $s_k(k_1)$ and $s_k(k_2)$, the first components are diffeomorphism types for smooth manifolds and the second components are both $1$.  
\end{enumerate}
For each of the two manifolds or polyhedra in (\ref{s:3c1}) and (\ref{s:3c2}), we have a canonical embedding into the newly obtained manifold or polyhedron in the PL category and we call such an embedding a {\it trace embedding}. Moreover, if the second component of $s_k(k_j)$ is $1$, then we say that the trace embedding is {\it special}.
\end{enumerate}
\item
\label{s:3}
Return to the first step replacing $k$ by $k+1$.
\end{enumerate}
\begin{Def}
\label{def:1}
We call a polyhedron in the first component of the resulting element $s_{l,0}$ and a polyhedron obtained in this way an {\it elementary polyhedron}.   
The sequence of the pairs of trace embeddings defined in (\ref{s:3c1}) or (\ref{s:3c2}) for such an iteration of these procedures is {\it associated with} the polyhedron. The length of the sequence is same as the time of procedures of defining the values as (\ref{s:3c1}) or (\ref{s:3c2}) in (\ref{s:3c}).

Last, we call the map $s_0$ and a map playing a same role in obtaining an elementary polyhedron a {\it root} for the polyhedron.
\end{Def}

\begin{Def}
\label{def:2}
Let $n \geq 1$ be an integer. An $n$-dimensional compact, connected and smooth manifold $X$ smoothly immersed (embedded) into ${\mathbb{R}}^n$ is said to be a {\it smoothly immersed }{\rm (}{\it embedded}{\rm )}{\it elementary manifold
 supported by $K$} or {\it SIE-$K$} if the following three are satisfied.
\begin{enumerate}
\item $X$ collapses to an elementary polyhedron $K$.
\item Any special embedding in a pair of a suitable sequence of trace embeddings associated with $K$ is smooth as the canonically obtained immersion into $X$. 
\item If for $K$ and a suitable procedure to obtain this before, the second component of $s_{l,0}$ is $1$, then the immersion of the smooth manifold $K$ into $X$ is also smooth.  
\end{enumerate}
\end{Def}
The $n$-dimensional compact manifold in Definition \ref{def:2} is regarded as a so-called {\it regular neighborhood} of $K$ in the smooth category. 

A systematic study of a regular neighborhood of a subolyhedron in a smooth manifold including the existence is discussed in \cite{hirsch} for example.
This notion implicitly appears in the proof of Main Theorem later for example.

As is also in \cite{kitazawa13}, construction in Theorem 6 of \cite{kitazawa8} can present various images of special generic maps explicitizing Definition \ref{def:2} well.
 
We end this section by a remark related to low-dimensional geometry which is closely related to our present study and which we do not concentrate on here. 
\begin{Rem}
In \cite{munozozawa} and \cite{ozawa} for example, {\it multibranched surfaces} are defined as $2$-dimensional polyhedra and these polyhedra embedded in $3$-dimensional compact manifolds are studied in the PL category. For example, such polyhedra represented as products of circles and graphs are studied. If the graph is obtained by taking bouquets starting from circles, then the multibranched surface is an elementary polyhedron. However, in the present paper, we do not consider such polyhedra since we concentrate on simply-connected polyhedra.
In $4$-dimensional (differential) topology, $4$-dimensional manifolds represented as so-called {\it handlebodies} and ones collapsing to $2$-dimensional polyhedra are important objects. For this, see \cite{gompgstipsicz} as an introductory book on $4$-dimensional (differential) topology for example.
\end{Rem}
\section{On Main Theorem.}

\begin{MainThm}
Let $n>k>1$ and $a>0$ be integers.
\begin{enumerate}
\item
\label{m:1}
Let $K$ be a ${\rm (}k-1{\rm )}$-connected polyhedron whose dimension is smaller than or equal to $n-a$.  
Assume that the values of a root $s_0$ for $K$ are always of the form $({[Y]}_{{\rm Diff}},1,p)$ satisfying either of the following two.
\begin{enumerate}
\item $({[Y]}_{{\rm Diff}},1,1)$ for a closed and connected manifold $Y$ which can be immersed {\rm (}resp. embedded{\rm )} smoothly into ${\mathbb{R}}^{\dim Y+a}$ with a trivial normal bundle.
\item $({[Y]}_{{\rm Diff}},1,0)$ for a closed and connected manifold $Y$ which can be immersed {\rm (}resp. embedded{\rm )} smoothly into ${\mathbb{R}}^{n}$.
\end{enumerate}

Then $Y$ here is always {\rm (}$k-1${\rm )}-connected, there exists an $n$-dimensional compact, connected and smooth manifold $X$ and $X$ is an SIE-$K$ {\rm (}resp. SEE-$K${\rm )}.
\item
\label{m:2}
For an $n$-dimensional compact, connected and smooth manifold $X$ and an elementary polyhedron $K$, suppose that $X$ is an SIE-$K$ {\rm (}SEE-$K${\rm )} and that the following two are satisfied.
\begin{enumerate}
\item $X$ and $K$ are {\rm (}$k-1${\rm )}-connected.
\item For a root $s_0$ for $K$, values are always of either of the following forms.
\begin{enumerate}
\item $({[Y]}_{{\rm Diff}},1,1)$ for a closed and connected manifold $Y$ which can be immersed {\rm (}resp. embedded{\rm )} smoothly into ${\mathbb{R}}^{\dim Y+a}$ with a trivial normal bundle.
\item $({[Y]}_{{\rm Diff}},1,0)$ for a closed and connected manifold $Y$ which can be immersed {\rm (}resp. embedded{\rm )} smoothly into ${\mathbb{R}}^{n}$.
\end{enumerate}
\end{enumerate}
Then $K$ is represented as a polyhedron obtained by a finite iteration of taking a suitable bouquet starting from a family of finitely many polyhedra each of which satisfies at least one of the following two.
\begin{enumerate}
\item A product of the following two polyhedra.
\begin{enumerate}
\item A closed and {\rm (}$k-1${\rm )}-connected manifold $F$ we can smoothly immerse {\rm (}resp. embed{\rm )} into ${\mathbb{R}}^{\dim F+a}$ with a trivial normal bundle.
\item A polyhedron $K_F$ such that for a root $s_{0,F}$ for $K_F$, values are always of the form $({[Y]}_{{\rm Diff}},1,1)$ for a closed and connected manifold $Y$ which can be immersed {\rm (}resp. embedded{\rm )} smoothly into ${\mathbb{R}}^{\dim Y+a}$ with a trivial normal bundle and that cannot be represented as a product of a closed and {\rm (}$k-1${\rm )}-connected manifold $Y^{\prime}$ we can smoothly immerse into ${\mathbb{R}}^{\dim Y^{\prime}+a}$ and another polyhedron $K_{Y,F}$.
\end{enumerate}
\item A closed and {\rm (}$k-1${\rm )}-connected manifold $F_0$ we can smoothly immerse {\rm (}resp. embed{\rm )} into ${\mathbb{R}}^{n}$.
\end{enumerate}
\end{enumerate}
\end{MainThm}
\begin{proof}[A proof of Main Theorem.]
We prove (\ref{m:1}). $Y$ is ($k-1$)-connected by virtue of the definitions of related notions such as an elementary polyhedron, an SIE-$K$ and an SEE-$K$.
To show the existence of a desired compact manifold $X$, it is sufficient to check the following two. They are fundamental facts and we give sketches of proofs of them.
\begin{itemize}
\item For a polyhedron $K$ obtained by a finite iteration of taking a suitable bouquet starting from $l \geq 1$ polyhedra $K_j$
for each of which we can take a ($\max\{\dim K_j\}+a$)-dimensional compact manifold regarded as an SIE-$K_j$ (resp. SEE-$K_j$) and taking a product of this and a closed and connected manifold $Y_0$
we can immerse (resp. embed) smoothly into ${\mathbb{R}}^{\dim Y_0+a}$ with a normal trivial bundle, there exists a ($\max\{\dim K_j\}+\dim Y_0+a$)-dimensional compact manifold which is an SIE-$K$ (resp. SEE-$K$).
\item For a manifold $K$ represented as a connected sum of closed, connected and orientable manifolds $Y_1$ and $Y_2$ we can immerse (resp. embed) smoothly into ${\mathbb{R}}^{\dim Y_j+a}$ with normal trivial bundles, there exists a ($\dim Y_j+a$)-dimensional compact manifold which is an SIE-$K$ (resp. SEE-$K$) and diffeomorphic to $K \times D^a$ for $j=1,2$. For a manifold $K$ represented as a connected sum of closed, connected and orientable manifolds $Y_1$ and $Y_2$ we can immerse (resp. embed) smoothly into ${\mathbb{R}}^{n}$, there exists an $n$-dimensional compact manifold which is an SIE-$K$ (resp. SEE-$K$).
\end{itemize}
First, for a closed, connected and orientable manifold $Z$ we can immerse (resp. embed) smoothly into the $a$-dimensional higher Euclidean space with a normal trivial bundle, there exists an $a$-dimensional higher compact manifold which is an SIE-$Z$ (resp. SEE-$Z$) and diffeomorphic to $Z \times D^{a}$. For the first property, a polyhedron $K^{\prime}$ obtained by a finite iteration of taking a suitable bouquet starting from $l$ polyhedra $K_j$
for each of which we can take a ($\max\{\dim K_j\}+a$)-dimensional compact manifold regarded as an SIE-$K_j$ (resp. SEE-$K_j$), we can easily obtain a ($\max\{\dim K_j\}+a$)-dimensional compact manifold which is an SIE-$K^{\prime}$ (resp. SEE-$K^{\prime}$). The ($\max\{\dim K_j\}+a$)-dimensional compact manifold which is an SIE-$K^{\prime}$ (resp. SEE-$K^{\prime}$) is regarded as a manifold represented as a boundary connected sum of the $l$ manifolds each of which is regarded as the ($\max\{\dim K_j\}+a$)-dimensional compact manifold regarded as an SIE-$K_j$ (resp. SEE-$K_j$). For the product $K=K^{\prime} \times Y_0$, we have a desired SIE-$K$ (resp. SEE-$K$), which is diffeomorphic to the product of $Y_0$ and the ($\max\{\dim K_j\}+a$)-dimensional compact manifold which is an SIE-$K^{\prime}$. For the second property, we can immerse $K$ (resp. embed) smoothly into ${\mathbb{R}}^{\dim Y_1+a}={\mathbb{R}}^{\dim Y_2+a}$ with a normal trivial bundle and we have a desired manifold, diffeomorphic to $K \times D^a$, in the former situation, and we have the similar fact in the latter situation.

We prove (\ref{m:2}). By virtue of the first assumption of the two assumptions, the fact that polyhedra we need to construct $K$ via the presented procedures are ($k-1$)-connected is shown in a similar way to the proof of the fact that $Y$ in (\ref{m:1}) is ($k-1$)-connected. $K$ is, by the definition, represented as a bouquet of polyhedra obtained via the presented procedures. For each polyhedron to obtain a desired polyhedron here, assume that this does not satisfy the latter property here. By virtue of the rule for the construction with the second assumption on a root, it is represented as a product of a closed and {\rm (}$k-1${\rm )}-connected manifold $F$ we can smoothly immerse into ${\mathbb{R}}^{\dim F+a}$ with a normal trivial bundle and another polyhedron. Last, the second assumption on a root also poses the restriction that $F_0$ can be smoothly immersed or embedded into ${\mathbb{R}}^n$ for the second property here. 

This completes the proof.
\end{proof}
Note that \cite{kitazawa13} concerns the case the following two are satisfied.
\begin{enumerate}
\item $n \leq 3k$.
\item The value of a root is always of the form $({[Y]}_{{\rm Diff}},1,p)$ where $Y$ is a homotopy sphere except for manifolds playing same roles as the manifold $''F_0''$ plays here.
\end{enumerate}
\begin{MainCor}
In Main Theorem, let $n=7$ and $k=2$. Then in {\rm (}\ref{m:2}{\rm )}, $F$, $F_0$ and $K_F$ satisfy the following properties.
\begin{enumerate}
\item $F$ is a $k_F$-dimensional closed and simply-connected manifold with $2 \leq k_F \leq 4$.
\item $K_F$ is represented as a bouquet of copies of unit spheres if $k_F \geq 3$ and in general represented as a bouquet of closed and simply-connected manifolds whose dimensions are smaller than $7-k_F$ and which we can smoothly immerse {\rm (}resp. embed{\rm )} into ${\mathbb{R}}^n$.
\item $F_0$ is represented as a connected sum of finitely many closed and simply-connected manifolds whose dimensions are smaller than $7$ we can smoothly immerse {\rm (}resp. embed{\rm )} into ${\mathbb{R}}^n$. 
\end{enumerate}
\end{MainCor}
Note that \cite{kitazawa13} concerns the case $n=5,6$. To know the topologies and the differentiable structures of $F$ and $F_0$ and the topology of $K_F$ more precisely, consult \cite{kitazawa13} and as differential topological theory
 of closed and simply-connected manifolds whose dimensions are at most $6$, \cite{barden}, \cite{gompgstipsicz}, \cite{milnor}, \cite{nishioka}, \cite{smale}, \cite{wall}, \cite{zhubr}, \cite{zhubr2}, and so on. As a related fact, every $4$-dimensional closed, connected and orientable manifolds can be smoothly embedded into ${\mathbb{R}}^7$.
\section{A conclusion remark.}
We have studied the topology and the differential structure of $W_f$, and the image of $W_f$ via $\bar{f}$, for a special generic map $f$ in Proposition \ref{prop:1}. 
Note again that this is a refined version of works in \cite{kitazawa13} for example. Proposition \ref{prop:1} seems to be a key ingredient in investigating the manifold $M$. This method works in \cite{nishioka} and contributes to proofs of some facts in Proposition \ref{prop:2}. 

Principles of this kind have been and will continue to be keys in investigating the world of closed (and simply-connected) manifolds via special generic maps, more generally, so-called {\it fold maps}, which are higher dimensional variants of Morse functions, and smooth maps of more general classes. 

Related to such forthcoming studies, we end the present paper by introducing \cite{kitazawa0.1}--\cite{kitazawa10}, which concern special generic maps and fold maps on explicit closed manifolds and are also introduced in the introduction and in \cite{kitazawa13}.

\section{Acknowledgement and on data.}
The author is a member of and supported by JSPS KAKENHI Grant Number JP17H06128 "Innovative research of geometric topology and singularities of differentiable mappings". The author declares that all data supporting the present study directly are all in this paper. 

\end{document}